\documentclass[a4paper]{article}

\usepackage{amsmath,amsthm,amsfonts,latexsym,amscd,amssymb,enumerate}
\usepackage[hypertex,backref]{hyperref} 

\newcommand{\N}{\mathbb{N}}
\newcommand{\T}{\mathbb{T}}
\newcommand{\Z}{\mathbb{Z}}
\newcommand{\Q}{\mathbb{Q}}

\newcommand{\map}{{\rm map}}
\newcommand{\res}{{\rm res}}

\theoremstyle{plain}
\newtheorem{theorem}{Theorem}
\newtheorem{corollary}[theorem]{Corollary}

\newtheorem{conjecture}[theorem]{Conjecture}
\newtheorem{example}[theorem]{Example}

\theoremstyle{definition}
\newtheorem{definition}[theorem]{Definition}

\theoremstyle{remark}

\newtheorem{rem}[theorem]{Remark}

\DeclareMathOperator{\im}{im}      

\begin{document}

\pagestyle{myheadings}
\markboth{Thomas Schick and Andreas Thom}{About a conjecture of D.~H.~Gottlieb}

\title{On a conjecture of Daniel H.\ Gottlieb}

\author{Thomas Schick\\ Mathematisches Institut\\
Bunsenstr. 3-5\\ 37073 G\"ottingen, Germany.
\and 
Andreas Thom\\ Mathematisches Institut\\
Bunsenstr. 3-5\\ 37073 G\"ottingen, Germany.}

\maketitle

\begin{abstract}
  We give a counterexample to a conjecture of D.H. Gottlieb and prove a
  strengthened version of it.

  The conjecture says that a map from a finite CW-complex $X$ to an
  aspherical CW-complex $Y$ with non-zero Euler characteristic can have
  non-trivial degree (suitably defined) only if the centralizer of the image
  of the fundamental group of $X$ is trivial.

  As a corollary we show that in the above
  situation all components of non-zero degree maps in the space of maps from
  $X$ to $Y$ are contractible.
\end{abstract}

\section{A version of Gottlieb's conjecture}

Let $X$ and $Y$ be finite CW-complexes. In \cite{go0,go1}, Gottlieb defines a notion of degree of a continuous
map $f\colon X \to Y$ as follows.
Let $f_*\colon {H_*}(X,\Z) \to {H_*}(Y,\Z)$ be the induced map in reduced integral homology.
The degree $\deg(f)$ of $f$ is the least integer $n \in \N$, such that there exists a group homomorphism
$\tau\colon {H_*}(Y,\Z) \to {H_*}(X,\Z)$ which satisfies $f_*\circ \tau = n \cdot {\rm id}$.
He conjectures the following (compare \cite{go1}):

\begin{conjecture}[Gottlieb] Let $(Y,y)$ be a finite aspherical CW-complex which is not acyclic. Let $f\colon (X,x) \to (Y,y)$ be a
continuous map with $\deg(f) \neq 0$. If $\chi(Y)\neq 0$, then the centralizer of $f_*(\pi_1(X,x))$ in $\pi_1(Y,y)$ is trivial.
\end{conjecture}

In this note we give a counterexample to this form of the conjecture (see
Example \ref{countex})  
and prove a version with a stronger hypothesis, see Theorem \ref{main}.
Let us rephrase one important consequence of \textit{non-vanishing degree} in the case of mappings between closed 
oriented manifolds,
so that it is applicable in a more general setting.

\begin{definition}
Let $f \colon (X,x) \to (Y,y)$ be a 
continuous map. We say that $f$ is a \emph{superposition},
if for any $\Q \pi_1(Y,y)$-module $L$, the induced map
$$f_* \colon H^{\pi_1(X,x)}_*(\tilde X; f^* L ) \to H^{\pi_1(Y,y)}_*(\tilde Y; L)$$
is surjective.
\end{definition}

We will see in Theorem \ref{exam}, that a map of non-vanishing degree between
closed oriented manifolds, or more generally between oriented Poincar\'e
duality complexes,
is a superposition. Moreover, an equivariant version of the Becker-Gottlieb transfer gives plenty of examples
of maps between CW-complexes which are not Poincar\' e
complexes. 

The problem with Gottlieb's definition of degree seems to be that it takes
only untwisted coefficients into account. Lead by Gottlieb, one can
therefore define a stronger version of degree as follows:
\begin{definition}
  A map $f\colon X\to Y$ between finite CW-complexes. Its \emph{twisted
    rational degree} $\deg_{tw,\Q}(f)$ is $1$ if $f$ is a superposition, and is
    $0$ otherwise.

    Its \emph{twisted degree} $\deg_{tw}(f)$ is the least positive integer
    $n\in\N$ such that for each $\Z\pi_1(Y)$-module $L$ there is
    a group homomorphism $\tau_L\colon H_*(Y,L)\to H_*(X,f^*L)$ such that
    $f_*\circ \tau_L=n\cdot \rm id$, or $0$ if no such integer exists.
\end{definition}

Clearly, a map of non-zero twisted degree is a superposition, so that the next
result shows that Gottlieb's conjecture is correct if one requires that the
\emph{twisted} degree is non-zero.\vspace{0.1cm} 

Our main result is the following.

\begin{theorem}\label{main}
Let $(Y,y)$ be a finite  aspherical CW-complex. Let $f\colon (X,x)
\to (Y,y)$ be a  
continuous superposition. If $\chi(Y)\neq 0$, then the centralizer of
$f_*(\pi_1(X,x))$ in $\pi_1(Y,y)$ is trivial.
\end{theorem}

Assuming Theorem \ref{main}, we can show some corollaries which generalize results from \cite{go1}.
Let $f$ be a continuous map from $X$ to $Y$. We denote by $\map(X,Y,f)$ the space of continuous maps
from $X$ to $Y$ which are homotopic to $f$.

\begin{corollary}
Let $Y$ be a finite aspherical CW-complex. Let $f\colon X \to Y$ be
a  
continuous superposition. If $\chi(Y)\neq 0$ and $Y$ is aspherical, then the 
mapping space $\map(X,Y,f)$ is contractible.
\end{corollary}
\begin{proof} 
If $Y$ is aspherical, then $\map(X,Y,f)$ is also aspherical because of the
following reasoning: we have to extend a given map from $S^n$ to $\map(X,Y,f)$
to $D^{n+1}$. By the exponential law, this means to extend a map from $X\times
S^n$ to $X\times D^{n+1}$. The latter space is obtained from the former by
attaching cells of dimension $n+1$ or higher. If $n\ge 2$, because
$\pi_k(Y)=0$ for $k\ge 2$, we can extend the map cell by cell as required.

Gottlieb showed in \cite{goz} 
that $\pi_1(\map(X,Y,f),f)$ is naturally isomorphic to the centralizer of $f_*(\pi_1(X,x))$ in $\pi_1(Y,y)$.
Hence the claim follows from Theorem \ref{main}.
\end{proof}
\begin{corollary}\label{corol:fis}
Let $(Y,y)$ be a finite aspherical CW-complex with $\chi(Y) \neq 0$. Every subgroup of finite index in $\pi_1(Y,y)$
has trivial centralizer. 
\end{corollary}
\begin{proof}
Let $G$ be a finite index subgroup of $\pi_1(Y,y)$. The induced map $f\colon BG \to Y$ is a superposition. 
Hence, Theorem \ref{main} implies the claim.
\end{proof}

\begin{proof}[Proof of Theorem \ref{main}]
Because $Y=B\pi_1(Y)$ is finite dimensional, $\pi_1(Y)$ is torsion
free. Therefore every non-trivial subgroup is infinite. Let us assume that the
centralizer of $f_*(\pi_1(X,x))$ in $\pi_1(Y,y)$ is infinite.

If $\pi$ is a discrete group, let $L\pi$ be its group von Neumann algebra. 
If $\chi(Y)$ is not zero, then the equivariant $L^2$-homology
$$H^{\pi_1(Y,y)}_*(\tilde Y; L \pi_1(Y,y))$$ cannot
be zero-dimensional in all degrees. Indeed,
$$0 \neq \chi(Y) = \sum_{k=0}^{\infty} (-1)^k \beta^{(2)}_k(Y),$$
by Atiyah's $L^2$-index theorem, see \cite[Theorem 6.80]{L}. Here $\beta^{(2)}_k(Y)$ denotes the 
$k$-th $L^2$-Betti number
$$\beta^{(2)}_k(Y) = \dim_{L \pi_1(Y,y)} H^{\pi_1(Y,y)}_k(\tilde Y; L \pi_1(Y,y)).$$

By assumption, the map $f\colon (X,x) \to (Y,y)$ is a superposition, so it induces a surjection
$$H^{\pi_1(X,x)}_k(\tilde X; f^*
L \pi_1(Y,y) ) \to H^{\pi_1(Y,y)}_k(\tilde Y; L \pi_1(Y,y) ),$$
for every $k \in \N$.
However, since $Y$ is aspherical, for every subgroup $G$ of $\pi_1(Y,y)$ which contains $f_{*}(\pi_1(X,x))$, this map 
can be factorized through
$$H^{G}_k(EG; \res^{\pi_1(Y,y)}_G L \pi_1(Y,y) ) = H^{G}_k(EG; LG)
\otimes_{LG} L \pi_1(Y,y).$$
Here, we used that $L\pi_1(Y,y)$ is flat as $LG$-module, compare \cite[Theorem
6.29]{L}. By the same theorem,
the $L \pi_1(Y,y)$-dimension of the right hand side is equal to $\beta^{(2)}_k(G)$. To derive a contradiction, it suffices to construct a
sub-group as above which has only vanishing $L^2$-Betti numbers.

\medskip
If the centralizer of $f_*(\pi_1(X,x))$ intersects non-trivially with $\pi_1(X,x)$, then the intersection is infinite,
since $\pi_1(Y,y)$ is torsion-free. In this case $f_*(\pi_1(X,x))$ has an infinite center and all its $L^2$-Betti numbers
are zero by Theorem $7.2$ of \cite{L}. If the intersection is trivial, we may pick a non-torsion element which centralizes
$f_*(\pi_1(X,x))$. Together with $f_*(\pi_1(X,x))$, it generates a copy of $f_*(\pi_1(X,x)) \times \Z$, which has
trivial $L^2$-Betti numbers by K\"unneth's Theorem, see Theorem $6.54\ (4)$ in \cite{L}. Hence, 
we arrive at a contradiction. \vspace{0.1cm}

Since $\pi_1(Y,y)$ is torsion-free, we conclude that the centralizer of $f_*(\pi_1(X,x))$ is trivial. This finishes the proof.
\end{proof}

\begin{rem}
  Note that, because the classifying space of $\im(f_*)\subset \pi_1(Y)$ is in
  general not a finite CW-complex, we have to use the generalization of
  $L^2$-Betti numbers to arbitrary CW-complexes of L\"uck as developed in
  \cite{MR1603853,MR1605818}. 
\end{rem}

The next theorem gives examples of maps which are superpositions.
\begin{theorem} \label{exam}
The following classes of maps are superpositions.
\begin{enumerate}
\item retractions,
\item continuous maps between oriented closed manifolds, or more generally Poincar\'e duality spaces, 
which have non-vanishing degree,
\item continuous maps $f\colon (X,x) \to (Y,y)$, whose homotopy fiber has the homotopy type
of a finite CW-complex and non-vanishing Euler characteristic.
\end{enumerate}
\end{theorem}
\begin{proof} The  statement about retractions follows from functoriality; for
  the second statement one uses the transfer given by Poincar\'e duality, and
  for the third the Gottlieb-Becker transfer (with twisted coefficients).
\end{proof}

\begin{rem}
  The proofs of the results presented so far show that the assumptions can be
  weakened as follows:
  \begin{itemize}
  \item $\chi(Y)\ne 0$ can be replaced by the assumption that $Y$ has at least
    one non-vanishing $L^2$-Betti number 
  \item The map $f$ being a superposition can be replaced by the assumption
    that $f$ induces a surjective homomorphism in $L^2$-cohomology.
  \end{itemize}
  That surjectivity in $L^2$-cohomology is true for inclusions of finite index
  subgroups and therefore Corollary
  \ref{corol:fis} holds under the weaker assumptions follows e.g.~from
  \cite{MR1828605}.
\end{rem}

\begin{rem}
  One should observe that Gottlieb's theorem, stating that the center of an
  aspherical finite CW-complex with non-trivial Euler characteristic is trivial, has been
  generalized considerably. Its strongest version now reads that such a group
  does not contain an infinite amenable normal subgroup.

  Our main application states that the centralizer of an image group is
  trivial; again we expect a generalization similar to the one about normal
  amenable subgroups. However, the correct notion of ``amenable centralizer''
  still has to be developed.
\end{rem}

\section{A counterexample to a strong form of the conjecture}
\label{sec:we-finish-this}

We finish this note by giving the desired counterexample to Gottlieb's Conjecture.
The tools in the construction are the techniques from the work of Baumslag, Dyer and
Heller, see \cite{BHD}.
\begin{theorem}[Baumslag-Dyer-Heller]  \label{asp}
There exists a finite aspherical and acyclic CW-complex $(D,\ast)$ whose fundamental group $\pi_1(D,\ast)$ contains
a copy of $\Z$.
\end{theorem}

\begin{example} \label{countex}
There exists a finite aspherical CW-complex $Y$ with $\chi(Y)=2$ and a continuous map $f \colon \T^2 \to Y$ which is
of degree one (taking Gottlieb's definition) and injective on the fundamental group. In particular, the centralizer
of $f_*(\pi_1(\T^2,\ast)) = \Z^2$ is infinite.
\end{example}
\begin{proof}[Construction of the Example] Let $D$ be as in Theorem \ref{asp}.
Starting with $\T^2$, we glue in two copies of $D$, 
along a circle in $D$ representing the generator of $\Z \subset \pi_1(D,\ast)$, and along each of the 
generators of the fundamental group of $\T^2$, to obtain the new space $Y$. Let $f\colon \T^2 \to Y$
be the natural map, induced by the glueing process.

This is exactly the type of construction of Baumslag-Dyer-Heller; since we 
glue along inclusions on the level of fundamental groups, the resulting 
space $Y$ is aspheric and 
the map $f_*\colon \Z^2 \to \pi_1(Y,y)$ is injective.
On the other hand, a look at the Mayer-Vietoris 
sequence shows that the map from $\T^2$ to $Y$ is an isomorphism 
in second integral homology, whereas $H_1(Y,\Z)=0$.
\end{proof}

\bibliographystyle{plain}
\bibliography{gottlieb}


\end{document}